\documentclass[12pt]{amsart}

\newcounter{point}
\usepackage{fp}
\usepackage{tikz}
\usepackage{xparse}
\usepackage{xstring}

\usetikzlibrary{positioning, calc, babel, snakes, patterns, cd}
\NewDocumentCommand{\findMax}{>{\SplitArgument{2}{ } }m}{\findMaxHelp#1}

\NewDocumentCommand{\findMaxHelp}{m m m} {
	\FPsub\h#2#1
	\FPadd\h\addconstant\h
	\FPmul\h\stretch\h
	\FPiflt\hcap\h\FPset\hcap\h\else\fi}

\FPset\lastheight{-3}
\FPset\drawingheight{0}
\FPset\currentmaxcapheight{0}
\FPset\stretch{0.5}
\FPset\addconstant{1}
\FPset\hcap{1}
\FPset\hcup{1}
\FPset\stddiff{1}
\FPset\scaling{0.5}
\NewDocumentCommand{\capbasic}{m m m}{
	\FPsub\h{#2}{#1}
	\FPadd\h\addconstant\h
	\FPmul\h\stretch\h
	\IfNoValueTF{#3}{\draw[-, black, out=270, in=270] (#1,{-\drawingheight}) .. controls+(0,{\h})  and +(0,{\h}) .. (#2, {-\drawingheight});}{
		\IfEq{#3}{d}{
			\draw[-, black, out=270, in=270] (#1,{-\drawingheight}) .. controls+(0,{\h})  and +(0,{\h}) .. (#2, {-\drawingheight})node[transform shape, circle, inner sep=3pt, midway, fill = black] {};}{\draw[-, black, out=270, in=270] (#1,{-\drawingheight}) .. controls+(0,{\h})  and +(0,{\h}) .. (#2, {-\drawingheight});}}
}

\NewDocumentCommand{\cupbasic}{m m m}{
	\FPsub\h{#2}{#1}
	\FPadd\h\addconstant\h
	\FPmul\h\stretch\h
	\IfNoValueTF{#3}{\draw[-, black, out=90, in=90] (#1,{-\drawingheight}) .. controls+(0,{-\h})  and +(0,{-\h}) .. (#2, {-\drawingheight});}{
		\IfEq{#3}{d}{
			\draw[-, black, out=90, in=90] (#1,{-\drawingheight}) .. controls+(0,{-\h})  and +(0,{-\h}) .. (#2, {-\drawingheight})node[transform shape, circle, inner sep=3pt, midway, fill = black] {} ;}{\draw[-, black, out=90, in=90] (#1,{-\drawingheight}) .. controls+(0,{-\h})  and +(0,{-\h}) .. (#2, {-\drawingheight});}}
	\FPiflt\hcup\h\FPset\hcup\h\else\fi
}

\NewDocumentCommand{\capintermediate}{>{\SplitArgument{2} { } }m }{\capbasic#1}

\NewDocumentCommand{\caps}{>{\SplitList{,}} m}{
	\FPset\hcap{1}
	\ProcessList{#1}{\findMax}
	\FPiflt\hcap\stddiff\FPset\diff\stddiff\else\FPset\diff\hcap\fi
	\FPifeq{-3}\lastheight{}\else\FPadd\drawingheight\diff\drawingheight\fi
	\ProcessList{#1}{\capintermediate}
}
\NewDocumentCommand{\cupintermediate}{>{\SplitArgument{2} { } }m }{\cupbasic#1}

\NewDocumentCommand{\cups}{>{\SplitList{,}} m}{
	\ProcessList{#1}{\cupintermediate}
	\FPset\lastheight\drawingheight
	\FPadd\drawingheight\hcup\drawingheight
}

\NewDocumentCommand{\raybasic}{m m m}{
	\IfNoValueTF{#3}{
		\raybasic{#1}{#2}{}
	}{
		\IfEq{#3}{d}{\IfEq{#1}{i}{
				\draw[-, black, out=90, in=270] (#2,{-\lastheight-\hcup}) -- (#2, {-\lastheight}) node[transform shape, midway, circle, inner sep = 3pt, fill=black]{};
			}{\IfEq{#2}{i}{
					\draw[-, black, out=90, in=270] (#1,{-\drawingheight}) -- (#1, {-\drawingheight+\hcap}) node[transform shape, midway, circle, inner sep = 3pt, fill=black]{};
				}{
					\draw[-, black, out=90, in=270] (#1,{-\drawingheight}) -- (#1, {-\drawingheight+\scaling*\hcap}) -- (#2, {-\lastheight-\scaling*\hcup}) node[transform shape, midway, circle, inner sep = 3pt, fill=black]{} -- (#2, {-\lastheight});
				}
		}}
		{
			\IfEq{#1}{i}{
				\draw[-, black, out=90, in=270] (#2,{-\lastheight-\hcup}) -- (#2, {-\lastheight});
			}{\IfEq{#2}{i}{
					\draw[-, black, out=90, in=270] (#1,{-\drawingheight}) -- (#1, {-\drawingheight+\hcap});
				}{
					\draw[-, black, out=90, in=270] (#1,{-\drawingheight}) -- (#1, {-\drawingheight+\scaling*\hcap}) -- (#2, {-\lastheight-\scaling*\hcup}) -- (#2, {-\lastheight});
				}
		}}
	}
}
\NewDocumentCommand{\rayintermediate}{>{\SplitArgument{2}{ }} m}{\raybasic#1}
\NewDocumentCommand{\rays}{>{\SplitList{,}} m}{
	\FPsub\diff\drawingheight\lastheight
	\FPiflt\diff\stddiff\FPsub\diff\stddiff\diff\FPadd\drawingheight\drawingheight\diff\else\fi
	\ProcessList{#1}{\rayintermediate}}

\NewDocumentCommand{\wdiag}{ >{\SplitList{ }} m }{
	\FPsub\diff\drawingheight\lastheight
	\FPiflt\diff\stddiff\FPsub\diff\stddiff\diff\FPadd\drawingheight\drawingheight\diff\else\fi
	\setcounter{point}{0}\ProcessList{#1}{\functions}\begin{pgfonlayer}{layer1}
		\draw (-0.5,{-\drawingheight})--(\value{point}-0.5, {-\drawingheight});
	\end{pgfonlayer}
	\FPset\lastheight\drawingheight
}

\NewDocumentCommand{\wdiagnoline}{ >{\SplitList{ }} m }{
	\FPsub\diff\drawingheight\lastheight
	\FPiflt\diff\stddiff\FPsub\diff\stddiff\diff\FPadd\drawingheight\drawingheight\diff\else\fi
	\setcounter{point}{0}\ProcessList{#1}{\functions}
	\FPset\lastheight\drawingheight
}

\NewDocumentCommand{\functions}{m}{
	\IfStrEqCase{#1}{
		{v}{\draw (\value{point}-.1, {-\drawingheight+.2}) -- (\value{point},{-\drawingheight}) -- (\value{point} +.1, {-\drawingheight+.2});}
		{w}{\draw (\value{point}-.1, {-\drawingheight-.2}) -- (\value{point},{-\drawingheight}) -- (\value{point} +.1, {-\drawingheight-.2});}
		{x}{\functions{v}\addtocounter{point}{-1}\functions{w}\addtocounter{point}{-1}}
		{o}{\node[circle, draw=black, fill=white, inner sep=-2pt, transform shape] at (\value{point}, {-\drawingheight}) {};}
		{d}{\begin{pgfonlayer}{layer0}
				\draw[fill=white] (\value{point}-.1, {-\drawingheight}) -- (\value{point},{-\drawingheight-.2}) -- (\value{point} +.1, {-\drawingheight})--(\value{point},{-\drawingheight+.2}) --  (\value{point}-.1, {-\drawingheight});
			\end{pgfonlayer}}
		{-}{}}[{\node[transform shape] at (\value{point}, {-\drawingheight}) {#1};}]
	\stepcounter{point}
}

\newcommand{\currh}{-\drawingheight}

\newcommand{\resetvariables}{
	\FPset\hcap{1}
	\FPset\hcup{1}
	\FPset\lastheight{-3}
	\FPset\drawingheight{0}
}

\tikzset{
	every scope/.append style={execute at begin scope={\resetvariables}},
}

\pgfdeclarelayer{layer0}
\pgfdeclarelayer{layer1}
\pgfsetlayers{main,layer0,layer1}
\usepackage[alphabetic, msc-links]{amsrefs}
\usepackage[utf8]{inputenc}
\usepackage{fontenc}
\usepackage{mathtools}
\usepackage{amssymb}
\usepackage{mathrsfs}
\usepackage{amsmath}
\usepackage{amsthm}
\usepackage{cleveref}
\usepackage{hyperref}
\usepackage{bm}
\usepackage{enumitem}
\usepackage{todonotes}
\usepackage{microtype}

\usepackage{caption}
\usepackage[labelfont={sf},font={small},
labelsep=space]{caption}
\calclayout

\crefname{equation}{}{}
\Crefname{equation}{Equation}{Equations}
\crefname{enumi}{}{}
\crefname{Def}{Definition}{Definitions}
\crefname{defi}{Definition}{Definitions}
\crefname{lem}{Lemma}{Lemmas}
\crefname{prop}{Proposition}{Propositions}
\crefname{thm}{Theorem}{Theorems}
\crefname{thmintro}{Theorem}{Theorems}
\crefname{cor}{Corollary}{Corollaries}
\crefname{rem}{Remark}{Remarks}
\crefname{section}{Section}{Sections}
\Crefname{section}{Section}{Sections}
\crefname{subsection}{Section}{Sections}
\Crefname{subsection}{Section}{Sections}
\crefname{figure}{Figure}{Figures}
\Crefname{figure}{Figure}{Figures}
\crefname{ex}{Example}{Examples}
\Crefname{ex}{Example}{Examples}

\newcommand{\bbZ}{\mathbb{Z}}

\newcommand{\bbC}{\mathbb{C}}
\DeclareMathOperator{\Mod}{mod}
\DeclareMathOperator{\gMod}{gmod}
\newcommand{\lie}[1]{\mathfrak{#1}}
\DeclareMathOperator{\Hom}{Hom}
\DeclareMathOperator{\mmod}{-mod}

\DeclareMathOperator{\Rep}{Rep}
\DeclareMathOperator{\cu}{cups}
\DeclareMathOperator{\B}{B}
\DeclareMathOperator{\A}{A}
\DeclareMathOperator{\ca}{caps}
\DeclareMathOperator{\Ext}{Ext}

\DeclareMathOperator{\id}{id}
\DeclareMathOperator{\im}{im}
\DeclareMathOperator{\Gl}{GL}
\DeclareMathOperator{\Or}{O}

\newcommand{\Osp}[1][m]{\def\ArgI{{#1}}\Osps}
\newcommand{\Osps}[1][n]{\mathrm{OSp}(\ArgI|#1)}
\newcommand{\KD}{\mathbb{D}}
\newcommand{\KB}{\mathbb{K}}

\newcommand{\cR}{\mathcal{R}}
\theoremstyle{definition}
\newtheorem{defi}{Definition}[section]
\newtheorem{rem}[defi]{Remark}

\newtheorem{ex}[defi]{Example}
\theoremstyle{plain}
\newtheorem{lem}[defi]{Lemma}

\newtheorem*{prop*}{Proposition}

\newtheorem{cor}[defi]{Corollary}
\newtheorem{thm}[defi]{Theorem}
\newtheorem*{thm*}{Theorem}

\DeclareMathOperator{\Br}{Br}
\DeclareMathOperator{\End}{End}

\begin{document}
\title{Koszulity for semi-infinite highest weight categories}

\author{{\rm Thorsten Heidersdorf and Jonas Nehme and Catharina Stroppel}}

\address{T. H.: School of Mathematics, Statistics and Physics, Newcastle University}
\email{heidersdorf.thorsten@gmail.com} 
\address{J. N.: Max-Planck-Institut für Mathematik, Bonn }
\email{nehme@mpim-bonn.mpg.de}
\address{C. S.: Mathematisches Institut, Universität Bonn}
\email{stroppel@math.uni-bonn.de}

\begin{abstract}
	We show that any upper finite or essentially finite highest weight category where the standard objects have linear projective resolutions and the costandard objects have linear injective resolutions is Koszul.
	This extends the result of \cite{ADL} to the case of infinite highest weight categories.
	We apply this result to Khovanov algebras and representations of classical Deligne categories and show that these are Koszul.
\end{abstract}

\maketitle 

\section{Introduction}

\subsection*{Koszul algebras} A Koszul algebra $A$ is a positively graded algebra with semisimple $A_0$ which is as ``close'' to being semisimple as possible, i.e.
\begin{equation*}
	\Ext^n_A(A_0, A_0\langle k\rangle) = 0\quad\text{unless}\quad n=k.
\end{equation*}  
This is equivalent to saying that all irreducible modules $L$ have a linear projective resolution, that is a projective resolution \begin{equation*}
	\begin{tikzcd}
		\dots\arrow[r]&P^2\arrow[r]&P^1\arrow[r]&P^0\arrow[two heads, r]&L
	\end{tikzcd}
\end{equation*}
such that $P^i$ is generated by its degree $i$ component (i.e.~$P^i=AP^i_i$, see \cite{BGS96} and \cite{MOS}). 

A Koszul algebra $A$ has a Koszul dual algebra $E(A)$. A main feature of the two is that their graded module categories are related, namely  (under additional finiteness assumptions)  the associated  graded derived categories of $A$ and of $E(A)$ are equivalent.  Koszulity is in general  a very desirable property, but often not easy to verify.   

In case the algebra is quasi-hereditary (which requires that it is finite dimensional), its categories of modules is a highest weight category in the sense of \cite{CPS}. In this situation, there is  a very convenient way to check Koszulity which is due to  \cite{ADL}: instead of showing the existence of linear projective resolutions for irreducible modules, it suffices to show it for standard modules. In practise this is usually much easier. One might even have extra tools available to construct these resolutions by induction on the partial ordering given by the highest weight structure.  A typical example is given by the algebra controlling a regular block of category $\mathcal{O}$ for a semisimple Lie algebra $\mathfrak{g}$, where translation functors allow such an induction argument. 

However he finiteness assumption required in these settings is very restrictive.

\subsection*{The main result} The main result of this article is an extension of the \cite{ADL} result to a semi-infinite setting. Motivated by the fact that many categories which exhibit highest weight structures satisfy natural, but weaker  finiteness conditions as in the classical setup, Brundan and Stroppel defined in \cite{BS21} three different types of highest weight categories called \emph{upper finite}, \emph{lower finite} and \emph{essentially finite}. 

A typical example is given by integral blocks of category $\mathcal{O}$ of positive, of negative and of critical level respectively for a Kac-Moody Lie algebra $\mathfrak{g}$.  By a result of Soergel, \cite{Soergel}, the case of positive and negative level are related by Ringel duality providing an example of the general result from \cite{BS21} that upper finite and lower finite highest weight categories are related by Ringel duality. As a consequence, it suffices to ask for a generalisation of the \cite{ADL} Koszulity criteria to the upper finite and to the essentially finite cases (the lower finite case follows then from the upper finite one by working with tilting resolutions instead of projective resolution). 

\begin{defi}
	Let $\cR$ be an upper finite or essentially finite highest weight category. We call $\cR$ \emph{standard Koszul} if every standard object $\Delta(\lambda)$ admits a linear projective resolution and every costandard object $\nabla(\lambda)$ admits a linear injective resolution.
\end{defi}

\begin{thm}\label{koszulcrit} 
	If $\cR$ is standard Koszul, then $\cR$ is Koszul.
\end{thm}

If $\cR$ admits a duality preserving the irreducibles, the condition on the $\nabla(\lambda)$ can be dropped. We omit the lower finite case here since these categories do not have enough projective objects unless they are already essentially finite.

\subsection*{Applications} We apply our result to two different categories and show that these are Koszul: the representation categories of Khovanov algebras of type A and B and of the Deligne categories $\Rep(\Or_{\delta})$ and $\Rep(\mathrm{GL}_{\delta})$, $\delta\in\mathbb{C}$.

Khovanov algebras were first considered in the context of categorification of tangle invariants \cite{K1}, \cite{Stroppel}. The generalized versions considered in this article were first introduced in \cite{BS1} for type $\A$ and in \cite{ES2} for type $\B$. These (generalized) Khovanov algebras have subsequently been used to introduce and study graded versions of blocks of $\mathfrak{gl}(m|n)$ \cite{BS4} and $\Osp[r][2n]$ \cite{ES3}, parabolic category $\mathcal{O}$ and the (walled) Brauer algebra \cite{BS5}.

Khovanov algebras depend on a set of weight diagrams, which can be finite or infinite in type $B$, $C$, $D$.
The finite dimensional versions in type $\A$ from \cite{BS1} and the one of type $\B$ from \cite{ES2} are known to be Koszul, which can be seen via identifying them with parabolic blocks of the Bernstein--Gelfand--Gelfand category $\mathcal{O}$ as in \cite{BS3}*{Theorem 1.1} (resp.~\cite{ES2}*{Theorem 9.1}), which are known to be Koszul (see e.g.~\cite{BGS96}). In \cite{BS2}, a different approach was used to directly prove that the Khovanov algebra of type $\A$ is Koszul without first identifying it with some category $\mathcal{O}$. We focus on the type $B$ version $\KB_{\Lambda}$ in the following. 

As modules over Khovanov algebras are semi-infinite highest weight categories, it suffices  by \cref{koszulcrit} to construct a linear projective resolution of the standard modules. The essential tool for this are projective functors, see also \cite{HNS24}. These functors are defined by tensoring with certain bimodules and correspond to translation functors when linking the Khovanov algebras to representations of $\Osp[r][2n]$ or to parabolic category $\mathcal{O}$. 

\begin{thm} For every standard module $\Delta(\lambda)$ there exists a linear projective resolution. In particular, $\KB_{\Lambda}$ is Koszul. 
\end{thm}

We refer to \cite{HNS24} for a detailed study of the involved modules.

The second set of examples arises from the theory of Deligne categories. 
Deligne categories are universal monoidal categories which interpolate categories of representations for families of groups like the classical groups $\Gl(m)$, $\Or(n)$ but also of supergroups like $\Gl(m|n)$ or $\Osp[r][2n]$ and allow to treat these cases as continuous families $\Rep(\mathrm{GL}_\delta)$ respectively $\Rep(\Or_\delta)$ with $\delta\in\bbC$. 
A representation of $\Rep(\Or_{\delta})$ is a contravariant functor from $\Rep(\Or_{\delta})$ to $Vect$, the category of finite dimensional complex vector spaces. 
We denote this category by $\mathcal{D}(\delta)$. By \cite[Corollary 2.11]{OSPII} the category $\mathcal{D}(\delta)$ is an upper finite highest weight category and there is an equivalence of categories between $\mathcal{D}(\delta)$ and $\gMod_{lf}(\KB_{\Lambda})$. 
Hence, we obtain:

\begin{cor}
	The graded version of the category of representations of the Deligne category $\Rep(\Or_\delta)$ is Koszul.
\end{cor}

Similar results hold in type $\A$.

\subsection*{Outlook} We expect that the results of this article can be applied in many more examples. One such example are path algebras of infinite quivers which allow truncations to finite quivers. A particularly beautiful example are the Cubist algebras by Chuang and Turner \cite{CT}. For these algebras Koszulity is already established in \cite{CT} by considering finite truncations of the underlying poset. (We need to retain minor finiteness assumptions as discussed in \cite{MOS} to exclude examples of quivers containing vertices with an infinite number of adjacent edges in which cases irreducible modules might not have any resolution with finitely generated projective modules at all.)

Other possible examples are representation categories of more general diagrammatic tensor categories since these often exhibit upper finite highest weight structures.

\subsection*{Acknowledgments}	
The authors are grateful to Hausdorff center for Mathematics in Bonn as well as Max Planck Institute for Mathematics in Bonn for its hospitality and financial support. The research of T.H., J.N. and C.S. was supported by the Deutsche For\-schungs\-ge\-mein\-schaft (DFG, German Research Foundation) under Germany's Excellence Strategy – EXC-2047/1 – 390685813.

\section{Koszul algebras and highest weight categories}
\begin{defi}
	A \emph{locally unital algebra $A=\bigoplus_{d\in\bbZ}A_d$} is an associative (but not necessarily unital) algebra together with a system $\{e_\lambda\mid\lambda\in\Lambda\}$ of mutually orthogonal idempotents such that 
	\begin{equation*}
		A=\bigoplus_{\lambda, \mu\in\Lambda}e_\lambda Ae_\mu.
	\end{equation*}
	
	The locally unital algebra $A$ is \emph{locally finite-dimensional} if $\dim e_\lambda A e_\mu<\infty$ for all $\lambda,\mu\in\Lambda$.
	It is called \emph{essentially finite-dimensional} if $\dim e_\lambda A$ and $\dim Ae_\lambda$ are finite for all $\lambda\in\Lambda$.
	By an $A$-module we always mean a left $A$-module $M$ such that 
	\begin{equation*}
		M=\bigoplus_{\lambda\in\Lambda}e_\lambda M.
	\end{equation*}
	An $A$-module $M$ is called \emph{locally finite-dimensional} if $\dim e_\lambda M<\infty$ for all $\lambda\in\Lambda$ and $d\in\bbZ$.
	If $A$ is locally unital, we will denote by $\Mod_{lf}(A)$ all locally finite dimensional $A$-modules and by $A\mmod$ the full subcategory of all finite dimensional ones.

\end{defi}
For the following definitions we will always assume that $A$ is a locally unital, locally finite dimensional positively graded algebra such that $A_0$ is semisimple.
Any module over $A$ is assumed to be locally finite dimensional.
\begin{defi}
	For $d\in\bbZ$, let $M$ be an $A$-module generated by its degree $d$ component.
	A linear projective resolution of $M$ is a projective resolution
	\begin{equation*}
		\begin{tikzcd}
			\dots\arrow[r]&P^2\arrow[r]&P^1\arrow[r]&P^0\arrow[two heads, r]&M
		\end{tikzcd}
	\end{equation*}
	such that $P^i$ is generated by its degree $i+d$ component (i.e.~$P^i=AP^i_{i+d}$).
	If this property only holds for $P^i$ with $i\leq n$ for some $n\in\bbZ_{\geq0}$, we call this an $n$-linear projective resolution.
	We also have the analogous notion of linear injective resolutions.
	The algebra $A$ is called \emph{Koszul} (see \cite{Priddy}, and also \cite{BGS96}, 
	\cite{MOS}) if every irreducible module admits a linear projective resolution.
\end{defi}	

We recall the definition of highest weight categories from \cite{BS21}.
\begin{defi}
	Let $A$ be a locally unital, locally finite dimensional algebra $A$ with semisimple $A_0$.
	Assume that $\cR$ is equivalent to either $\Mod_{lf}(A)$ or $A\mmod$ (if $A$ is essentially finite-dimensional). 
	A poset $\Lambda$ is called \emph{upper finite} if $[\mu,\infty)$ is finite for all $\mu \in \Lambda$, \emph{lower finite} if $(\infty,\mu]$ is finite for all $\mu\in\Lambda$ and \emph{essentially finite} if $[\mu,\lambda]$ is finite for all $\lambda,\mu \in \Lambda$.
	Suppose that $\Lambda$ is upper finite in case $\cR \cong \Mod_{lf}(A)$ and essentially finite in case $\cR \cong A\mmod$, and suppose that there is a bijection $\lambda\mapsto L(\lambda)$ from $\Lambda$ onto the isomorphism classes of simple objects up to grading shift.
	
	We have Serre subcategories $\cR_{\leq \lambda}$ and $\cR_{<\lambda}$ of $\cR$ generated by all simple objects $L(\mu)$ with $\mu\leq\lambda$ and $\mu<\lambda$, respectively.
	We say that $\cR$ is an \emph{upper finite} (respectively \emph{essentially finite}) \emph{highest weight category} if the following conditions are satisfied:
	\begin{enumerate}
		\item The Serre quotient category $\cR_\lambda\coloneqq\cR_{\leq\lambda}/\cR_{<\lambda}$ is simple for all $\lambda\in\Lambda$.
		\item Denote by $j^{\lambda}\colon\cR_{\leq\lambda} \to \cR_\lambda\coloneqq\cR_{\leq\lambda}/\cR_{<\lambda}$ the quotient functor and set $P_\lambda(\lambda)=L_{\lambda}(\lambda) = j^{\lambda}L(\lambda)$ and denote by $P_{\lambda}(\lambda)$ its projective cover. 
		The functor $j^{\lambda}$ has a left adjoint $j^{\lambda}_{!}$ and we set $\Delta(\lambda) = j_{!}^{\lambda} P_\lambda(\lambda)$. 
		Then we require: For every $\lambda\in\Lambda$ the projective cover $P(\lambda)$ of $L(\lambda)$ has a finite filtration with $\Delta(\lambda)$ at the top and other sections $\Delta(\mu)$ for $\mu\geq\lambda$.
	\end{enumerate}
\end{defi}

\begin{ex} The following categories provide examples of semi-infinite highest weight categories.
\begin{enumerate}
\item \cite[Theorem 6.4]{BS21} Let $\mathfrak{g}$ be an affine Kac--Moody algebra and let $\mathcal{O}_\Lambda$ be an integral block of $\mathcal{O}$ of non-critical level. In non-critical levels, the irreducible modules in the block are parametrized by the elements in the Weyl group orbit $W\cdot\lambda$ where $\lambda$ is integral of non-critical level. In case the level is positive, the block $\mathcal{O}_{\Lambda}$ is an upper finite highest weight category, and the standard objects are the Verma modules $\Delta(\lambda)$ for $\lambda \in \Lambda$. The partial order $\leq$ on $\Lambda$ is the dominance order.
\item Various types of Deligne or interpolation categories give rise to upper finite highest weight categories. While these categories are not necessarily abelian, one can consider their category of representations, contravariant functors from $\mathcal{C}$ to $Vect$ (finite dimensional vector spaces). An example is given by the partition category studied in \cite{D19}, \cite{BV}. This category of representations is an upper finite highest weight category. A similar situation occurs in the case of the Deligne category $\Rep(\Or_{\delta})$, $\delta \in \mathbb{C}$ (see Section \ref{sec4} for details). It was also shown in \cite{GRS} that one gets upper finite highest weight categories in this way from the cyclotomic Brauer category and the cyclotomic Kauffman category.
\end{enumerate}
\end{ex}
\begin{rem}
	To obtain a \emph{graded highest weight category}, the following changes have to be made:
	\begin{enumerate}
		\item We consider graded locally unital algebras $A$ where the pieces $e_\lambda A e_\mu$ are degree-wise finite dimensional and bounded below.
		All modules are assumed to be graded and locally finite dimensional, meaning that $e_\lambda M_d$ is finite dimensional for all $\lambda\in\Lambda$ and $d\in\bbZ$.
		\item The filtration of the projective cover $P(\lambda)$ is not necessarily finite anymore, but we require that $\Delta(\mu)\langle k\rangle$ appears in the filtration only for finitely many $\mu\in\Lambda$.
		Furthermore, the multiplicity of $\Delta(\mu)\langle k\rangle$ is finite and $0$ for $k\ll 0$.
	\end{enumerate}
\end{rem}
\begin{defi}\label{standardKoszul}
	Let $\cR$ be an upper finite or essentially finite graded highest weight category. 
	We call $\cR$ \emph{standard Koszul} if every standard object $\Delta(\lambda)$ admits a linear projective resolution and every costandard object $\nabla(\lambda)$ admits a linear injective resolution.
\end{defi}
\begin{rem}\label{standardKoszulwithduality}
	If $\cR$ admits a duality preserving the irreducibles, then the existence of linear injective resolutions of costandard objects follows from applying the duality to the linear projective resolutions of standard objects.
	So in this case, it suffices to check the existence of linear projective resolutions of standard objects.
\end{rem}
\begin{rem}
	We excluded the case of $\cR$ being a lower finite highest weight category, as these usually do not have projective objects.
	In fact, if a lower finite highest weight category has projective objects, it is actually already essentially finite.
\end{rem}
The main theorem of this paper is the following theorem, extending the result of \cite{ADL} to the case of infinite highest weight categories.
\begin{thm}\label{koszulcrit}
	If $\cR$ is standard Koszul, then $\cR$ is Koszul.
\end{thm}
Before we are going to prove this, we will introduce some preparatory results.

\begin{lem}[Homological criterion for linear projective resolution]\label{homcrit}
	Consider $\Delta(\nu)$ for $\nu\in\Lambda$.
	Then $\Delta(\nu)$ admits a linear projective resolution if and only if $\Ext^r_A(\Delta(\nu), L(\mu)\langle k\rangle)=0$ unless $r=d-k$ for all $\mu\in\Lambda$.
\end{lem}
\begin{proof}
	Consider a projective resolution of $\Delta(\nu)$:
	\begin{equation*}
		\begin{tikzcd}
			\dots\arrow[r]&P^2\arrow[r]&P^1\arrow[r]&P^0\arrow[two heads, r]&\Delta(\nu).
		\end{tikzcd}
	\end{equation*}
	We have $P^i=\bigoplus_{\lambda, k}(P(\lambda)\langle k\rangle)^{\oplus n_{\lambda, k}}$ for some $n_{\lambda, k}\in\bbZ_{\geq0}$.
	From the assumptions on the $\Delta$-filtration of any indecomposable projective, we may assume that only finitely many $P(\lambda)$ appear with degree shift $k$ (and none appear with degree shift $k$ for $k\ll 0$).
	
	Since any degree $0$ morphism $f\colon P(\lambda)\to P(\mu)$ implies that $f$ is an isomorphism (as $A_0$ semisimple), we can inductively remove these terms (starting from $i=0$ and small $k$).
	In particular, we can assume that all differentials have degree $>0$.
	Therefore, if we apply $\Hom_A(\_, L(\lambda)\langle k\rangle)$ to the complex, every differential vanishes.
	Thus, we have $\Ext^r_A(\Delta(\nu), L(\lambda)\langle k\rangle)=0$ unless $r=d-k$ if and only if $\Delta(\nu)$ has a linear projective resolution.
\end{proof}


This homological criterion allows us to show that standard Koszulity translates to upper truncations and lower truncations.
We assume in the following that $\cR$ is an upper finite or essentially finite highest weight category with poset $\Lambda$ (which is upper finite respectively interval finite).

\begin{lem}\label{stdkoszulrestrictionupper}
	Let $\Lambda^{\uparrow}\subseteq\Lambda$ be an upper set (i.e.~$\lambda\in\Lambda^{\uparrow}$ implies $\mu\in\Lambda^{\uparrow}$ for all $\mu\geq\lambda$).
	Let $\cR^{\uparrow}$ be the corresponding Serre quotient category (i.e.~the quotient by the Serre subcategory generated by $L(\lambda)$ for $\lambda\notin\Lambda^{\uparrow}$).	
	If $\cR$ is standard Koszul, then $\cR^{\uparrow}$ is standard Koszul.
\end{lem}
\begin{proof}
	Write $j\colon\cR\to\cR^{\uparrow}$ for the quotient functor and denote by $j_!$ the left adjoint functor and by $j_*$ the right adjoint.
	By \cite{BS21}*{Thm.~3.18} and \cite{BS21}*{Thm.~3.42} respectively, we have the following properties
	\begin{itemize}
		\item $\cR^{\uparrow}$ is an upper finite or essentially finite highest weight category with simple objects $L^{\uparrow}(\lambda)\coloneqq j L(\lambda)$, standard objects $\Delta^{\uparrow}(\lambda)\coloneqq j\Delta(\lambda)$ and costandard objects $\nabla^{\uparrow}(\lambda)\coloneqq j\nabla(\lambda)$ for $\lambda\in\Lambda^{\uparrow}$,
		\item $j_!\Delta^{\uparrow}(\lambda) = \Delta(\lambda)$ for all $\lambda\in\Lambda^{\uparrow}$,
		\item $\Ext^r(j_!\Delta^{\uparrow}(\lambda), L(\mu)\langle k\rangle)=\Ext^r(\Delta^{\uparrow}(\lambda), j L(\mu)\langle k\rangle)$ for all $r\geq 0, k\in\bbZ$ and $\lambda$, $\mu\in\Lambda$.
		\item $j_*\nabla^{\uparrow}(\lambda) = \nabla(\lambda)$ for all $\lambda\in\Lambda^{\uparrow}$,
		\item $\Ext^r(L(\mu)\langle k\rangle, j_*\nabla^{\uparrow}(\lambda),)=\Ext^r(j L(\mu)\langle k\rangle, \nabla^{\uparrow}(\lambda))$ for all $r\geq 0, k\in\bbZ$ and $\lambda$, $\mu\in\Lambda$.
	\end{itemize}
	In virtue of \cref{homcrit}, we have to show that $\Ext^r(\Delta^{\uparrow}(\lambda), L^{\uparrow}(\mu)\langle k\rangle)=0$ unless $r=k$ for all $\lambda, \mu\in\Lambda^{\uparrow}$.
	We have  
	\begin{equation*}
		\Ext^r(\Delta^{\uparrow}(\lambda), L^{\uparrow}(\mu)\langle k\rangle)=\Ext^r(j_!\Delta^{\uparrow}(\lambda), L(\mu)\langle k\rangle)=\Ext^r(\Delta(\lambda), L(\mu)\langle k\rangle)
	\end{equation*}
	which is zero unless $r=k$ by assumption.
	The same argument shows the statement for the linear injective resolutions of costandard objects (using the analogue of \cref{homcrit}).
\end{proof}

\begin{lem}\label{stdkoszulrestrictionlower}
	Let $\Lambda^{\downarrow}\subseteq\Lambda$ be a lower set (i.e.~$\lambda\in\Lambda^{\downarrow}$ implies $\mu\in\Lambda^{\downarrow}$ for all $\mu\leq\lambda$).
	Let $\cR^{\downarrow}$ be the corresponding Serre subcategory.
	If $\cR$ is standard Koszul, then $\cR^{\downarrow}$ is standard Koszul.
\end{lem}

\begin{proof}
	Write $i\colon\cR^{\downarrow}\to\cR$ for the inclusion functor and denote by $i^*$ the left adjoint functor and by $i^!$ the right adjoint.
	By \cite{BS21}*{Thm.~3.17} and \cite{BS21}*{Thm.~3.41} respectively, we have the following properties
	\begin{itemize}
		\item $\cR^{\downarrow}$ is an upper finite or essentially finite highest weight category with simple objects $L^{\downarrow}(\lambda)\cong L(\lambda)$, standard objects $\Delta^{\downarrow}(\lambda)\cong\Delta(\lambda)$ and costandard objects $\nabla^{\downarrow}(\lambda)\cong\nabla(\lambda)$ for $\lambda\in\Lambda^{\downarrow}$,
		\item $i^*\Delta(\lambda) = \Delta^{\downarrow}(\lambda)$ for all $\lambda\in\Lambda^{\downarrow}$,
		\item $\Ext^r(i^*\Delta(\lambda), L(\mu)\langle k\rangle)=\Ext^r(\Delta(\lambda), i L(\mu)\langle k\rangle)$ for all $r\geq 0, k\in\bbZ$ and $\lambda,\mu\in\Lambda$.
		\item $i^!\nabla(\lambda) = \nabla^{\downarrow}(\lambda)$ for all $\lambda\in\Lambda^{\downarrow}$,
		\item $\Ext^r(L(\mu)\langle k\rangle, i^!\nabla(\lambda))=\Ext^r(i L(\mu)\langle k\rangle, \nabla^{\downarrow}(\lambda))$ for all $r\geq 0, k\in\bbZ$ and $\lambda,\mu\in\Lambda$.
	\end{itemize}
	In virtue of \cref{homcrit}, we have to show that $\Ext^r(\Delta^{\downarrow}(\lambda), L^{\downarrow}(\mu)\langle k\rangle)=0$ unless $r=k$ for all $\lambda, \mu\in\Lambda^{\downarrow}$.
	We have  
	\begin{equation*}
		\Ext^r(\Delta^{\downarrow}(\lambda), L^{\downarrow}(\mu)\langle k\rangle)=\Ext^r(i^*\Delta(\lambda), L^\downarrow(\mu)\langle k\rangle)=\Ext^r(\Delta(\lambda), L(\mu)\langle k\rangle)
	\end{equation*}
	which is zero unless $r=k$ by assumption.
	The same argument shows the statement for the linear injective resolutions of costandard objects.
\end{proof}

\subsection{Upwards induction}
In this section we show the main theorem in the case of an essentially finite highest weight category where the partial order is additionally lower finite instead of interval finite.
We start from a simple stratum, where the statement is trivial, and then move upwards by extending the partially ordered labeling set with bigger elements.
This approach is inspired by \cite{ADL} and generalizes the arguments there.
However, we will need to adapt the arguments as we cannot assume finite global dimension.

So assume throughout that $\cR$ is a standard Koszul essentially finite highest weight category with lower finite labeling set.
We need the following technical definition.
\begin{defi}
	Let $\lambda\in\Lambda$ be minimal (which exists as $\Lambda$ is lower finite).
	Denote by $j\colon \cR\to\cR_{\nleq\lambda}$ the quotient functor and by $j_!$ the left adjoint functor.
	We denote by $\kappa^\lambda_n$ the objects $X$ of $\cR$ satisfying the following properties:
	\begin{itemize}
		\item $X$ is generated in degree $d$.
		\item $X'\coloneqq \im(j_!jX\to X)$ is also generated in degree $d$.
		\item $X$ is $L(\lambda)$-Koszul, i.e.~$\Ext^r(X, L(\lambda)\langle k\rangle)=0$ unless $r=k-d$.
		\item $jX$ has an $n$-linear projective resolution, i.e.~a projective resolution that is linear in the first $n$ terms.
	\end{itemize}
\end{defi}
Next we show that any object in $\kappa^\lambda_n$ admits an $n$-linear projective resolution.

\begin{lem}\label{kappalinproj}
	Every module in $\kappa^\lambda_n$ has an $n$-linear projective resolution.
\end{lem}
\begin{proof}
	We prove the statement via induction on $n$.
	For $n=0$ the statement is trivial.
	Now suppose that the statement is true for $<n$.
	
	We have to show that $X$ admits an $n$-linear projective resolution, and we distinguish two cases.
	Either $X'\neq X$ or $X'=X$.
	We reduce the first case to the second case and the second case to our induction hypothesis.
	
	So assume that $X'\neq X$.
	By the triangle identities of adjunctions, we have $jX'=jX$.
	Therefore, as $j$ is exact, $j(X/X') = 0$ and thus $X/X'$ is of the form $\bigoplus L(\lambda)\langle k\rangle$.
	As $X$ is generated in degree $d$, we have $X/X'$ generated in degree $d$ as well, so $X/X'$ is of the form $\bigoplus L(\lambda)\langle d\rangle$.
	Then, we have a short exact sequence
	\begin{equation}\label{ses}
		\begin{tikzcd}[column sep = small]
			0\arrow[r]&X'\arrow[r]&X\arrow[r]&\bigoplus L(\lambda)\langle d\rangle\arrow[r]&0.
		\end{tikzcd}
	\end{equation}
	As $\lambda$ was chosen minimal, we have $L(\lambda)=\Delta(\lambda)$.
	Hence, as $\cR$ is standard Koszul, $L(\lambda)$ has an $n$-linear projective resolution.
	By the horseshoe lemma, $X$ has an $n$-linear projective resolution if $X'$ has an $n$-linear projective resolution (as $X'$ is also generated in degree $d$).
	As $jX'=jX$, we get $X''=X'$ and if we can show that $X'\in\kappa^\lambda_n$ it suffices to show the statement for $Y'=Y$.
	Using $jX'=jX$, we get that $jX'$ has an $n$-linear projective resolution and clearly $X''$ is generated in degree $d$.
	So we only need to show $X'$ is $L(\lambda)$-Koszul.
	The short exact sequence \cref{ses} induces a long exact sequence
	\begin{align*}
		\dots\to&\Ext^r(\bigoplus L(\lambda)\langle d\rangle, L(\lambda)\langle k\rangle)\to\Ext^r(X, L(\lambda)\langle k\rangle)\to\Ext^r(X', L(\lambda)\langle k\rangle)\to\\&\Ext^{r+1}(\bigoplus L(\lambda)\langle d\rangle, L(\lambda)\langle k\rangle)\to\dots.
	\end{align*}	
	As $\lambda$ minimal, we have $\Delta(\lambda)\cong\nabla(\lambda)\cong L(\lambda)$ and hence $\Ext^r( L(\lambda)\langle d\rangle, L(\lambda)\langle k\rangle)=0$ for all $r\geq 1$ by \cite{BS21}*{Cor.~3.2}.
	Together with the $L(\lambda)$-Koszulity of $X$ we get that $X'$ is $L(\lambda)$-Koszul.
	Hence, $X'\in\kappa$.
	
	The remaining case to consider is $X'=X$.
	Consider the short exact sequence
	\begin{equation*}
		\begin{tikzcd}[column sep = small]
			0\arrow[r]&\Omega\arrow[r]&P\arrow[r]&X\arrow[r]&0
		\end{tikzcd},
	\end{equation*} where $P$ is the projective cover of $X$ and $\Omega$ the first syzygy.
	If we can show that $\Omega$ is generated in degree $d+1$ and that $\Omega\in\kappa^\lambda_{n-1}$, we can apply the induction hypothesis to get that $\Omega$ has an $(n-1)$-linear projective resolution.
	Then, we can assemble this with the above short exact sequence into an $n$-linear projective resolution of $X$.
	So we have to show that $\Omega\in\kappa^\lambda_{n-1}$ and that it is generated in degree $d+1$.
	As $j$ is exact, we also obtain an exact sequence
	\begin{equation*}
		\begin{tikzcd}[column sep = small]
			0\arrow[r]&j\Omega\arrow[r]&jP\arrow[r]&jX\arrow[r]&0
		\end{tikzcd}.
	\end{equation*}
	By assumption $jX$ has an $n$-linear projective resolution and thus $j\Omega$ has an $(n-1)$-linear projective resolution.
	In particular, we see that $j\Omega$ is generated in degree $d+1$ and thus also $\Omega'$.
	Now consider the projective cover $P'$ of $\Omega$.
	We have to show that $P'$ is generated in degree $d+1$.
	For every summand of the form $P(\mu)\langle k\rangle$ with $\mu\neq\lambda$ this follows from the statement for $\Omega'$.
	For $P(\lambda)\langle k\rangle$ this follows from the $L(\lambda)$-Koszulity of $X$.
	Finally, we have to show that $\Omega$ is $L(\lambda)$-Koszul.
	
	For this consider the long exact sequence
	\begin{align*}
		\dots\to&\Ext^r(P, L(\lambda)\langle k\rangle)\to\Ext^r(\Omega, L(\lambda)\langle k\rangle)\to\Ext^{r+1}(X, L(\lambda)\langle k\rangle)\to\\&\Ext^{r+1}(P, L(\lambda)\langle k\rangle)\to\dots.
	\end{align*}
	The outer terms vanish for all $r>1$ as $P$ is projective.
	As $X=X'$, no summand of $P$ can be isomorphic to $P(\lambda)\langle k\rangle$.
	Hence, $\Ext^r(\Omega, L(\lambda)\langle k\rangle)\cong\Ext^{r+1}(X, L(\lambda)\langle k\rangle)$, which is $0$ unless $r+1=k-d$.
	Thus, $\Omega$ is $L(\lambda)$-Koszul.
	Finally, this gives $\Omega\in\kappa^\lambda_{n-1}$.
	By induction hypothesis $\Omega$ has an $(n-1)$-linear projective resolution and thus we can assemble this with the above short exact sequence into an $n$-linear projective resolution of $X$.
\end{proof}

\begin{rem} In Lemma \ref{kappalinproj} we are using that we work with a highest weight category (not just a stratified category). Indeed, we need the existence of a linear projective resolution of an arbitrary object in $\mathcal{R}_{\lambda}$ for minimal $\lambda$. If $\mathcal{R}_{\lambda}$ is simple, the existence is clear.
\end{rem}

Now we have all the ingredients to prove the main theorem for the special case of an essential finite highest weight category $\cR$ with lower finite labeling set.
\begin{thm}\label{koszulcritlowerfinite}
	Let $\cR$ be an essentially finite highest weight category with lower finite labeling set.
	If $\cR$ is standard Koszul, then $\cR$ is Koszul.
\end{thm}
\begin{proof}
	We need to show that every $L(\mu)$ has a linear projective resolution.
	Let $\lambda\leq\mu$ minimal.
	If $\mu=\lambda$, we have $L(\lambda)=\Delta(\lambda)$ and hence a linear projective resolution by assumption.
	If $\mu\neq\lambda$, we want to reduce the statement to $\cR_{\nleq\lambda}$.
	For this, observe that $L(\mu)$ as well as $L(\mu)'$ are generated in degree $0$.
	
	Next, we consider $\Ext^r(L(\mu), L(\lambda)\langle k\rangle)$.
	By assumption, we know that $L(\lambda)=\nabla(\lambda)$ has a linear injective resolution.
	But this means, that $\Ext^r(L(\mu), L(\lambda)\langle k\rangle)=0$ unless $r=k$ and hence $L(\mu)$ is $L(\lambda)$-Koszul.
	Hence, in view of \cref{kappalinproj}, $L(\mu)$ has a linear projective resolution if $j L(\mu)$ has a linear projective resolution, where $j\colon\cR\to\cR_{\nleq\lambda}$ is the quotient functor.
	
	By \cref{stdkoszulrestrictionupper}, $\cR_{\nleq\lambda}$ is standard Koszul with lower finite labeling set.
	Furthermore, $j L(\mu)$ is again an irreducible module in $\cR_{\nleq\lambda}$.
	So we can repeat the argument from above.
	As $\Lambda$ is lower finite, we will reach a step where $\mu$ is minimal, whence we get a linear projective resolution.
	As mentioned above this gives then linear projective resolutions for the previous Serre quotient categories and ultimately for $L(\mu)$.
\end{proof}
\subsection{Downwards induction}
In this section, we prove the main theorem.
The main strategy is to use \cref{stdkoszulrestrictionupper} to restrict to essentially finite highest weight categories with lower finite labeling set.
In the restricted case, we can apply \cref{koszulcritlowerfinite} to show that $\cR^{\uparrow}$ is Koszul.

Then, we assemble all the upper truncations into a directed system for $\cR$ and show that the colimit gives the desired linear projective resolution.
\begin{thm}
	Let $\cR$ be an upper finite or essentially finite graded highest weight category, respectively.
	
	If $\cR$ is standard Koszul, then $\cR$ is Koszul.
\end{thm}
\begin{proof}
	We have to show that every $L(\lambda)$ has a linear projective resolution.
	For this, we first introduce some technical upper subsets of $\Lambda$.
	Let $X\subseteq\Lambda$ and write $B_1(X)$ for the set of all $\lambda\in\Lambda$ such that there exists a $\mu\in X$ with $\Hom_{\cR}(P(\lambda)\langle1\rangle, P(\mu))\neq 0$.
	Define $B_n(X)\coloneqq B_{n-1}(B_1(X))$ and let $\Lambda_n$ be the upper set generated by $B_n(\{\lambda\})$.
	As $\Hom_{\cR}(P(\lambda)\langle1\rangle, P(\mu))\neq 0$ is nonzero for only finitely $\lambda$ (for fixed $\mu$) by assumption, we, in particular, have that $\Lambda_n$ is lower finite.

	We write $\cR_n$ for the upper truncation of $\cR$ at $\Lambda_n$.
	For $\lambda\in\Lambda_n$ we write ${e_nL(\lambda)}$ for the irreducible associated to $\lambda$ in $\cR_n$.
	We use similar notation for the standard and projective modules.
	
	By $j^n\colon \cR\to \cR_n$ we denote the Serre quotient functor.
	This has a left adjoint $j^n_!\colon \cR_n\to\cR$ and a right adjoint $j^{n,*}\colon\cR_n\to \cR$ given by $\Hom_{A_n}(e_nA, \_)$.
	
	Furthermore, for $n\leq m$ we have Serre quotient functors $j_m^n\colon \cR_m\to\cR_n$ and the corresponding left respectively right adjoints $j_{m,!}^n$ and $j_m^{n,*}$.
	The left adjoint is called \emph{standardization functor} and the right adjoint is called \emph{costandardization functor}
	These satisfy the following commutation relations ($n<m<k$):
	\begin{equation}\label{commrelations}
		\begin{aligned}
			j_m^n\circ j_k^m&=j_k^n&		j_m^n\circ j^m&=j^n\\		
			j^m_{k,!}\circ j_{m,!}^n&=j^n_{k,!}&		j^m_!\circ j_{m,!}^n&=j^n_!\\
			j^{m,*}_k\circ j^{n,*}_m&=j^{n,*}_k&		j^{m,*}\circ j^{n,*}_m&=j^{n,*}\\
		\end{aligned}
	\end{equation}
	Additionally for $\lambda\in\Lambda_m$ and $n<m$ the functor $j^n_m$ sends ${e_mP(\lambda)}$ respectively ${e_mL(\lambda)}$ to ${e_nP(\lambda)}$ respectively ${e_nL(\lambda)}$ if $\lambda\in\Lambda_n$ and ${e_mL(\lambda)}$ to $0$ otherwise.
	The same holds if we leave out the $m$.
	
	The standardization functor $j_{m,!}^n$ sends ${e_nP(\lambda)}$ to ${e_mP(\lambda)}$ for $\lambda\in\Lambda_n$, and we have that $j_m^n\circ j_{m,!}^n=\id$.
	Again the same holds true if we leave out the $m$.

	By assumption $\cR$ is standard Koszul and hence $\cR_n$ is standard Koszul for all $n$ by \cref{stdkoszulrestrictionupper}.
	The Serre quotient $\cR_n$ is an (essentially) finite highest weight category with (lower) finite labeling set.
	Hence, it is Koszul by \cref{koszulcritlowerfinite}.
	So for $\lambda\in\Lambda_m$ we have a linear projective resolution of ${e_mL(\lambda)}$ in $\cR_m$, which is unique up to isomorphism.
	We denote this by (where $P_m^k(\lambda)$ is generated in degree $k$)
	\begin{equation*}
		P^\bullet_m(\lambda)\colon \dots\to P_m^k(\lambda)\to P_m^{k-1}(\lambda)\to\dots\to P_m^1(\lambda)\to P_m^0(\lambda)\to {e_mL(\lambda)}.
	\end{equation*}

	Now if $\lambda\in\Lambda_n$ for some $n<m$ by exactness of $j^n_m$ we get an exact sequence
	\begin{equation*}
		\dots\to j^n_mP_m^k(\lambda)\to j^n_mP_m^{k-1}(\lambda)\to\dots\to j^n_mP_m^1(\lambda)\to j^n_mP_m^0(\lambda)\to {e_nL(\lambda)}.
	\end{equation*}
	Therefore, by \cite{Wei}*{Thm.~2.2.6} there is a map (unique up to homotopy) \begin{equation*}
		\iota^n_m(\lambda)\colon P^{\bullet}_n(\lambda)\to j^n_mP^{\bullet}_m(\lambda),
	\end{equation*} and we can choose these maps $\iota^n_m(\lambda)$ such that for $n\leq m\leq k$ the following diagram commutes
	\begin{center}
		\begin{tikzcd}
			P^\bullet_n(\lambda)\arrow[r, "\iota^n_m(\lambda)"]\arrow[dr, "\iota^n_k(\lambda)"]&j^n_mP_m^\bullet(\lambda)\arrow[d, "j^n_m(\iota^m_k(\lambda))"]\\
			&j^n_kP^\bullet_k(\lambda)=j^n_mj^m_kP^\bullet_k(\lambda)
		\end{tikzcd}.
	\end{center}
	For $\mu\in\Lambda_m$ and $\mu\notin\Lambda_n$ for $m>n$ we have that ${e_mP(\mu)}$ can only appear in homological degrees $>n$ of a linear projective resolution of ${e_mL(\lambda)}$ by construction of $\Lambda_n$.
	Therefore, any summand of $P_m^k(\lambda)$ for $k\leq n$ has to be of the form ${e_mP(\mu)}$ for some $\mu\in\Lambda_n$.

	Additionally, $j^n_m{e_mP(\mu)}={e_nP(\mu)}$ for $\mu\in\Lambda_n$ and thus,
	\begin{center}
		\begin{tikzcd}
			j^n_mP^{n}_m(\lambda)\arrow[r]&\dots\arrow[r]&j^n_mP^{1}_m(\lambda)\arrow[r]&j^n_mP^{0}_m(\lambda)\arrow[r]&{e_nL(\lambda)}
		\end{tikzcd}
	\end{center}
	is a beginning of a linear projective resolution of ${e_nL(\lambda)}$ and by uniqueness of a linear projective resolution, $\iota^n_m(\lambda)$ has to be an isomorphism in degrees $\leq n$.
	
	Using the unit $\eta^n_m$ of the adjunction $j_{m,!}^n\vdash j_m^n$, we get a morphism $g_m^n(\lambda)$ of resolutions as the composition
	\begin{center}
		\begin{tikzcd}
			j_{m,!}^nP^\bullet_n(\lambda)\arrow[rr, "j_{m,!}^n(\iota_m^n(\lambda))"]&&j_{m,!}^nj_m^nP^\bullet_m(\lambda)\arrow[rr, "(\eta^n_m)_{P^\bullet_m(\lambda)}"]&& P^\bullet_m(\lambda).
		\end{tikzcd}
	\end{center}
	We observe that for $\mu\in\Lambda_n$ we have $j_{m,!}^nj_m^n{e_mP(\mu)}={e_mP(\mu)}$, so $g_m^n(\lambda)$ is an isomorphism in degrees $\leq n$.
	Additionally, we look at the diagram
	\begin{center}
		\begin{tikzcd}
			j^n_{k,!}P^\bullet_n(\lambda)=j^m_{k,!}j^n_{m,!}P^\bullet_n(\lambda)\arrow[r, "j^n_{k, !}(\iota^n_m(\lambda))"]\arrow[d, "j^n_{k,!}(\iota^n_k(\lambda))"]&j^m_{k,!}j^n_{m,!}j^n_mP^\bullet_m(\lambda)\arrow[r, "j^m_{k,!}(\eta^n_m)"]\arrow[d, "j^n_{k,!}j^n_m(\iota^m_k(\lambda))"]&j^m_{k,!}P^\bullet_m(\lambda)\arrow[d, "j^m_{k,!}(\iota^m_k(\lambda))"]\\
			j^n_{k,!}j^n_kP^\bullet_k(\lambda)\arrow[r, "="]\arrow[dr, "\eta^n_k"]&j^m_{k,!}j^n_{m,!}j^n_mj^m_kP^\bullet_k(\lambda)\arrow[r, "j^m_{k, !}(\eta^n_m)"]&j^m_{k,!}j^m_kP^\bullet_k(\lambda)\arrow[dl, "\eta^m_k"]\\
			&P^\bullet_k(\lambda)&
		\end{tikzcd}.
	\end{center}
	In order to make the diagram a bit more clear we suppressed for the unit $\eta^n_m$ the index describing the object.
	By our choice of the maps $\iota^n_m(\lambda)$ the left square commutes.
	By naturality of $\eta^n_m$ the right square commutes as well.
	The lower triangle commutes as the adjunction $j^n_{k,!}\vdash j^n_k$ is given as the composition of the adjunctions $j^n_{m,!}\vdash j^n_m$ and $j^m_{k,!}\vdash j^m_k$.
	Now there are two possibilities to go from $j^n_{k,!}P^\bullet_n(\lambda)$ to $P^\bullet_k(\lambda)$.
	One is given by $g^n_k(\lambda)$ and the other one by $g^m_k(\lambda)\circ j^m_{k,!}g^n_m(\lambda)$.
	So we obtain $g^n_k(\lambda)=g^m_k(\lambda)\circ j^m_{k,!}g^n_m(\lambda)$.
	
	By defining $f_m^n(\lambda)\coloneqq j_!^mg_m^n(\lambda)$ we get morphisms of chain complexes in $\cR$.
	\begin{center}
		\begin{tikzcd}
			j^n_!P^\bullet_n(\lambda)\arrow[r, "f_m^n(\lambda)"]&j^m_!P^\bullet_m(\lambda).
		\end{tikzcd}
	\end{center}
	
	From $g^n_k(\lambda)=g^m_k(\lambda)\circ j^m_{k,!}g^n_m(\lambda)$ we get $f_m^n\circ f_k^m=f_k^n$  for $n\leq m\leq k$.
	
	We claim that the direct limit of this sequence is a linear projective resolution of $L(\lambda)$.
	First note that the morphism $f_m^n(\lambda)$ is an isomorphism in homological degrees $\leq n$ (because $g^n_m(\lambda)$ is one in these degrees).
	This means that in homological degrees $k\leq n$ the direct limit is actually a projective module which is generated in degree $k$.
	So in each homological degree, we have by construction a finite direct sum of indecomposable projectives, and so this is a chain complex of locally finite dimensional modules.
	
	Let us take a closer look at this resulting chain complex
	\begin{equation*}
		\dots\to \varinjlim j^n_!P^s_n(\lambda)\to\dots\to\varinjlim j^n_!P^{1}_n(\lambda)\to\varinjlim j^n_!P^0_n(\lambda)\to\varinjlim j^n_!L(\lambda)\to 0.
	\end{equation*}
	An object $M$ in $\cR$ is zero if and only if $j^mM=0$ for all $m\geq 0$.
	As taking homology commutes with exact functors, we observe that the above chain complex is exact if and only if $j^m\varinjlim j^n_!P^\bullet_n$ is exact for all $m\geq 0$.
	
	The restriction functor $j^m$ admits a left adjoint ($j^{m,*}$) and so it commutes with direct limits.
	Hence, we have 
	\begin{equation*}
		j^m\varinjlim j^n_!P^s_n(\lambda)=\varinjlim j^mj^n_!P^s_n(\lambda).
	\end{equation*}
	Using the commutations relations \eqref{commrelations}, we see that \begin{equation*}
		j^mj^n_!=\begin{cases}
			j^n_{m,!}&\text{if $n<m$,}\\
			\id&\text{if $n=m$ and}\\
			j^m_n&\text{if $n>m$,}
		\end{cases}
	\end{equation*} and thus the right-hand side of the above equation is the direct limit of the complex
	\begin{equation*}
		j^0_{m,!}P^s_0(\lambda)\to\dots\to j^{m-1}_{m,!}P^s_{m-1}\to P^s_m(\lambda)\to j^m_{m+1}P^s_{m+1}(\lambda)\to j^m_{m+2}P^s_{m+2}(\lambda)\to\dots.
	\end{equation*}
	Now by our previous argument, the homomorphisms $j^m_{l}P^{s}_l(\lambda)\to j^m_{l+1}P^{s}_{l+1}(\lambda)$ are isomorphisms for $s\leq m$.
	Thus, the first $m$ terms in $j^m\varinjlim j^n_!P^\bullet_n(\lambda)$ agree with the first $m$ terms in $j^m_{l}P^\bullet_l(\lambda)$ and thus are exact.
	Hence, $\varinjlim j^n_!P^\bullet_n(\lambda)$ is also exact, which is then a linear projective resolution of $L(\lambda)$.

	Hence, $\cR$ is Koszul.
\end{proof}

In the following two sections we are going to apply our main theorem to show the Koszulity of two upper finite highest weight categories: modules over Khovanov algebras of type $B$ and graded representation categories of Deligne categories.

\section{Definition of the Khovanov algebra}\label{sec3}

In this section we are going to recall the definition of the Khovanov algebra of type $\B$ from \cite{ES2} and its most important properties. 
For the easier algebra in type $\A$ we refer to \cite{BS1}.
For this we fix $\delta\in\bbZ$ and $L_\infty\coloneqq\frac{\delta}{2}+\bbZ\cap\mathbb{R}_{\geq 0}$.
\begin{defi}
	An (infinite) weight diagram of type $\B$ is a map $\mu\colon L_\infty\to\{\wedge,\vee,\circ,\times,\diamond\}$ satisfying
	\begin{itemize}
		\item $\diamond$ can only occur as the image of $0$,
		\item $0$ can be only mapped to $\circ$ or $\diamond$,
		\item $\mu^{-1}(\{\circ, \times, \wedge\})$ is finite,
		\item $|\mu^{-1}(\circ)|-|\mu^{-1}(\times)|=\lfloor\frac{\delta}{2}\rfloor$.
	\end{itemize}
	
	Two weight diagrams belong to the same block, if the underlying number lines, the positions of $\circ$ and $\times$ and the number of $\wedge$'s agree.
\end{defi}
\begin{defi}Two vertices of a weight diagram are \emph{neighbored} if they are only separated by $\circ$ and $\times$. 
	The cup diagram $\underline{\mu}$ associated with a weight diagram $\mu$ is obtained by applying the following steps in order.
	\begin{enumerate}[label=(C-\arabic*)]
		\item First connect neighbored vertices labeled $\vee\wedge$ successively by a cup, i.e.~we connect the vertices by an arc forming a cup below.
		Repeat this step as long as possible, ignoring already joint vertices.
		Note that the result is independent of the order in which the connections are made
		\item Attach a vertical ray to each remaining $\vee$.
		\item\label{enum:dotone} Connect pairs of neighbored $\wedge$'s from left to right by cups (we interpret $\diamond$ for this as a $\wedge$).
		\item\label{enum:dottwo} If a single $\wedge$ or $\diamond$ remains, attach a vertical ray.
		\item\label{enum:dotthree} Put a marker $\bullet$ on each cup created in \cref{enum:dotone} and each ray created in \cref{enum:dottwo}.
		\item We erase the marker from the component that contains the $\diamond$ if the number of placed markers in \cref{enum:dotthree} is odd.
		\item Finally, delete all $\vee$ and $\wedge$ labels at vertices.
	\end{enumerate} 
	
	A cap diagram is defined to be the horizontal mirror image of a cup diagram.
	If $a$ is a cup or a cap diagram, its horizontal mirror image is denoted by $*$.
	We will also write $\overline{\lambda}$ for $\underline{\lambda}^*$.
\end{defi}
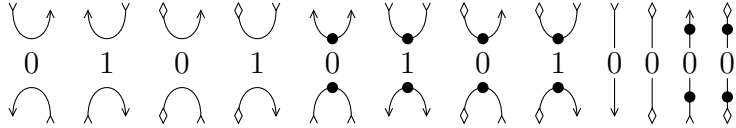
\begin{figure}[h]
	\centering
	\begin{tikzpicture}[scale=0.5]
		\FPset\stddiff{2}
		\wdiagnoline{v w w v d w d v w w v v d w d v v d w d}\cups{0 1, 2 3, 4 5, 6 7, 8 9 d, 10 11 d, 12 13 d, 14 15 d}\rays{i 16, i 17, i 18 d, i 19 d}\node at (0.5, {\currh+0.6}) {0};\node at (2.5, {\currh+0.6}) {1};\node at (4.5, {\currh+0.6}) {0};\node at (6.5, {\currh+0.6}) {1};
		\node at (8.5, {\currh+0.6}) {0};\node at (10.5, {\currh+0.6}) {1};\node at (12.5, {\currh+0.6}) {0};\node at (14.5, {\currh+0.6}) {1};
		\node at (16, {\currh+0.6}) {0};\node at (17, {\currh+0.6}) {0};\node at (18, {\currh+0.6}) {0};\node at (19, {\currh+0.6}) {0};
		\begin{scope}[yshift={2*\currh cm + 1.2cm}]
			\caps{0 1, 2 3, 4 5, 6 7, 8 9 d, 10 11 d, 12 13 d, 14 15 d}\rays{16 i, 17 i, 18 i d, 19 i d}\wdiagnoline{v w w v d w d v w w v v d w d v v d w d}
		\end{scope}
		\FPset\stddiff{1}
	\end{tikzpicture}
	\caption{Orientations and degrees}\label{orient}
\end{figure}
\begin{defi}
	An oriented circle diagram $a\lambda b$, is given by a cup diagram $a$, a weight diagram $\lambda$ and a cap diagram $b$ with the same underlying number lines, such that the positions of $\circ$ and $\times$ in $a$, $\lambda$ and $b$ agree and that each cup (resp.~cap) in $a$ (resp.~$b$) is oriented as in \cref{orient}.
	
	Every cup, cap and ray of $a\lambda b$ has an associated integer according to \cref{orient}.
	The sum $\deg(a\lambda b)$ of all these integers is called the \emph{degree} of $a\lambda b$.
\end{defi}
\begin{defi}\label{defkhovanovalg}
	For a block $\Lambda$ associated to the number line $L_\infty$, the Khovanov algebra $\KB_\Lambda$ is the graded associative algebra with underlying basis given by all oriented circle diagrams $a\lambda b$ with $\lambda\in\Lambda$ with underlying number line $L$ (resp. $L_\infty$), where $a\lambda b$ is homogeneous of degree $\deg(a\lambda b)$.
	
	The multiplication $(a\lambda b)(c\mu d)$ is defined to be $0$ whenever $b^*\neq c$ and if $b^*=c$, we draw the circle diagram $(a\lambda b)$ under the circle diagram $(b^*\mu d)$, where we connect the rays of $b$ and $b^*$ and apply certain surgery procedures.
	All these surgery procedures take a cup-cap-pair and replace it by two straight lines.
	After every cup-cap-pair is removed, one collapses the middle section and defines this to be $(a\lambda b)(c\mu d)$.
	For details (and the proof that this is in fact associative) we refer to \cite{ES2}*{Section 5}.
\end{defi}

\begin{rem}
	There also exists a finite dimensional version of this Khovanov algebra, see e.g.~\cite{ES2}.
	All properties and results mentioned here for $\KB_\Lambda$ hold for the finite dimensional one as well (using the same proofs).
\end{rem}
%
%
%
%
The algebra $\KB_\Lambda$ is not finite dimensional.
However, the elements $e_{\lambda}\coloneqq\underline{\lambda}\lambda\overline{\lambda}$ for $\lambda\in\Lambda$ form again pairwise orthogonal idempotents and thus give the algebra $\KB_\Lambda=\bigoplus_{\lambda, \mu\in\Lambda}e_\lambda \KB_\Lambda e_\mu$ the structure of a locally unital algebra.
By $\gMod_{lf}(\KB_\Lambda)$ we refer to locally finite dimensional graded modules over $\KB_\Lambda$, i.e.~graded modules $M$ such that $\dim e_\lambda M<\infty$ for all $\lambda\in\Lambda$.

The irreducible locally finite dimensional $\KB_\Lambda$-modules are in bijection with $\Lambda$.
Given $\lambda\in\Lambda$ we construct a one dimensional irreducible $\KB_\Lambda$-module $L(\lambda)$ as follows.
As a vector space it is just $\bbC$ and $e_\mu$ acts by $1$ if $\lambda=\mu$ and $0$ otherwise.

The indecomposable projective objects in $\gMod_{lf}(\KB_\Lambda)$ are given by $P(\lambda)\coloneqq \KB_\Lambda e_\lambda$ for $\lambda\in\Lambda$.

Furthermore, we define \emph{standard modules} $\Delta(\mu)$ for $\mu\in\Lambda$.
These are the cell modules associated to the cellular structure (in the sense of \cite{GL}) of $\KB_\Lambda$ in \cite{ES2}*{Theorem 7.1}.
These are given as the quotient of $P(\mu)$ by the submodules $U(\mu)$ generated by all oriented circle diagrams $a\lambda\overline{\mu}$ with $\lambda\neq\mu$ (then we necessarily have $\lambda>\mu$).

Given a graded $\KB_\Lambda$-module $M$, we define for $i\in\bbZ$ the graded module $M\langle i\rangle$ by $M\langle i\rangle_j = M_{j-i}$.

We have an anti-involution $*$ on $\KB_\Lambda$ which is given by sending $a\lambda b$ to $b^*\lambda a^*$.
And this gives rise to a duality (also denoted $*$) on $\gMod_{lf}(\KB_\Lambda)$.
For a locally finite dimensional graded $\KB_\Lambda$-module $M$, we define the graded piece $(M^{\circledast})_j\coloneqq\Hom_{\bbC}(M_{-j}, \bbC)$ and $x\in \KB_\Lambda$ acts on $f\in M^{\circledast}$ by $(xf)(m)\coloneqq f(x^*m)$.
We also easily see that $L(\lambda)^{\circledast}=L(\lambda)$.

With these definitions we can conclude this section with the following theorem:
\begin{thm}\label{upperfinitehighweight}
	The category $\gMod_{lf}(\KB_\Lambda)$ is an upper finite highest weight category in the sense of \cite{BS21} with standard objects $\Delta(\lambda)$, $\lambda\in\Lambda$.
\end{thm}
\begin{proof}
	The second part is just \cite{OSPII}*{Cor 2.11} after identifying their category $\mathcal{D}(\delta)$ (which consists of contravariant functors from $\Rep(\Or_\delta)$ to $\mathrm{Vect}$, the category of finite dimensional complex vector spaces) with $\gMod_{lf}(K)$ using \cite{OSPII}*{Thm.~6.22}.
\end{proof}
\begin{rem}
	The main difference between the type $\A$ and type $\B$ versions is that for type $\A$ there do not appear any dots on the circle diagrams but all the other properties listed here hold similarly in this easier case (see also \cite{BS1}).
\end{rem}

\section{Koszulity of the Khovanov algebras of type \texorpdfstring{$\B$}{B}}\label{sec4}
In this section we will show that $\KB_\Lambda$ is Koszul for any block $\Lambda$.
Note that in \cite{ES2}*{Theorem 9.1}, Ehrig and Stroppel identified the finite dimensional version with parabolic blocks of category $\mathcal{O}$, which are known to be Koszul (see e.g.~\cite{BGS96}).
This approach does not carry over to the infinite dimensional case, so we have to apply other methods.

As $\KB_\Lambda$ has a admits a duality that preserves the simple objects, it suffices to construct a linear projective resolution of the standard modules by \cref{koszulcrit}.
We will construct these using geometric bimodules and projective functors similar to the approach of Brundan and Stroppel in \cite{BS2}.
The definition of geometric bimodules uses crossingless matchings (of type $\B$) which we recall now.

A \emph{crossingless matching} is a diagram $t$, which is obtained by drawing an admissible cap diagram $c$ underneath an admissible cup diagram $d$ and connecting the rays in $c$ to the rays in $d$ from left to right.
This means that we allow dotted cups, caps and lines, but each dot necessarily needs to be able to be connected to the left boundary without crossing anything, just as in the case of admissible circle diagrams (see \cite{ES2}*{Def. 3.5}).
Furthermore, we delete pairs of dots on each segment, such that each line segment contains at most one dot.
Any crossingless matching is a union of (dotted) cups, caps and line segments, for example:
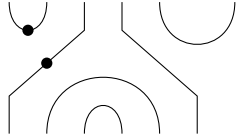
\begin{center}
	\begin{tikzpicture}[scale=0.5]
		\cups{1 2 d, 5 7}
		\caps{2 5, 3 4}
		\rays{1 3 d, 6 4}
	\end{tikzpicture}
\end{center}
We denote by $\cu(t)$ respectively $\ca(t)$ the number of cups respectively caps in $t$.
Furthermore, let $t^*$ be the horizontally reflected image of $t$.

We say that $t$ is a \emph{$\Lambda\Gamma$-matching} if the bottom and top number lines of $t$ agree with the number lines of $\Lambda$ respectively $\Gamma$.

Given additionally a cup diagram $a$ and a cap diagram $b$ such that their number lines agree with the bottom respectively top number line of $t$, we can glue them together and obtain a \emph{$\Lambda\Gamma$-circle diagram} $atb$.

Let $\Lambda$ and $\Gamma$ be blocks and let $t$ be a $\Lambda\Gamma$-matching.
Given weights $\lambda\in\Lambda$ and $\mu\in\Gamma$ we can glue these together from bottom to top to obtain a new diagram $\lambda t\mu$.
We call this an \emph{oriented $\Lambda\Gamma$-matching} if
\begin{itemize}
	\item each pair of vertices lying on the same dotted cup or the same undotted line segment is labeled such that both are either $\vee$ or both are $\wedge$,
	\item each pair of vertices lying on the same undotted cup or the same dotted line segment is labeled such that one is $\vee$ and one is $\wedge$,
	\item all other vertices are labeled $\circ$ or $\times$.
\end{itemize}
A diamond $\diamond$ can be interpreted as either $\vee$ or $\wedge$.

Finally, given an oriented $\Lambda\Gamma$-matching and cap and cup diagrams $a$ and $b$ such that $a\lambda$ (resp.~$\mu b$) is an oriented cup (resp cap) diagram we can glue these together to obtain an \emph{oriented $\Lambda\Gamma$-circle diagram} $a\lambda t\mu b$.

We call a $\Lambda\Gamma$-matching $t$ \emph{proper} if there exists at least one oriented $\Lambda\Gamma$-matching for $t$.
By a rightmost vertex $x$ on a circle $C$ we mean a vertex lying on $C$ such that on this numberline, there is no vertex to the right of $x$.
In the bottom picture every rightmost vertex is marked by $x$.
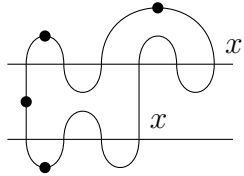
\begin{center}
	\begin{tikzpicture}[scale=0.5]
		\caps{0 1 d, 2 5 d, 3 4}\wdiag{- - - - - -}\node[anchor=south west] at (5, \currh) {$x$};\cups{1 2, 4 5}\caps{1 2}\rays{0 0 d, 3 3}\wdiag{- - - - - -}\node[anchor=south west] at (3, \currh) {$x$};\cups{0 1 d, 2 3}
	\end{tikzpicture}
\end{center}
We refer to a circle in an oriented $\Lambda\Gamma$-diagram as \emph{clockwise} respectively \emph{anticlockwise} if a rightmost vertex on the circle is labeled $\vee$ respectively $\wedge$.
This notion is well-defined by a similar argument as in \cite{ES2}*{ Corollary 5.9}.

\begin{defi}
	The \emph{degree} of a circle or a line in an oriented $\Lambda\Gamma$-circle diagram is the total number of clockwise cups or caps that it contains.
	The \emph{degree} of an oriented $\Lambda\Gamma$-circle diagram is the sum of the degrees of each of its circles and lines.
	We call a circle only consisting of one cup and one cap a \emph{small circle}.
\end{defi}
\begin{lem}\label{onemoreoneless}
	The degree of an anticlockwise circle in an oriented $\Lambda\Gamma$-circle diagram is one less than the total number of caps (equivalently, cups) that it contains.
	The degree of a clockwise circle is one more than the total number of caps (equivalently, cups) that it contains.
	The degree of a line is equal to the number of caps or the number of cups that it contains, whichever is greater.
\end{lem}
\begin{proof}
	This can be verified similar to \cite{ES2}*{Prop. 1.2.12+1.2.13}
\end{proof}

\begin{defi}\label{defiuplowreduction}
	Let $t$ be a $\Lambda\Gamma$-matching for some blocks $\Lambda$ and $\Gamma$.
	Let $a$ be a cap diagram such that its number line agrees with the top one of $t$.
	We refer to circles or lines not meeting the bottom number line in $ta$ as \emph{upper} circles or lines.
	The \emph{upper reduction} of $ta$ refers to the cup diagram which is obtained by removing all upper circles and lines as well as the top number line.
\end{defi}
\begin{ex}
	Suppose $ta$ as in \cref{defiuplowreduction} looks like:
	
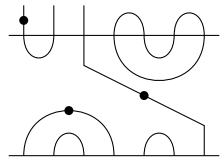
\begin{center}
		\begin{tikzpicture}[scale=0.4, yscale=1]
			\caps{3 4, 5 6}\rays{0 i d, 1 i, 2 i}\wdiag{- - - - - - -}\cups{0 1, 3 6, 4 5}\caps{0 3 d, 1 2, 4 5}\rays{6 2 d}\wdiag{- - - - - - -}
		\end{tikzpicture}
	\end{center}
	Then the upper reduction is:
	
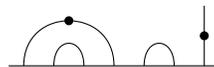
\begin{center}
		\begin{tikzpicture}[scale=0.4, yscale=1]
			\caps{0 3 d, 1 2, 4 5}\rays{6 i d}\wdiag{- - - - - - -}
		\end{tikzpicture}
	\end{center}
\end{ex}
\begin{lem}\label{degreeformulaupperreduction}
	If $a\lambda t\mu b$ is an oriented $\Lambda\Gamma$-circle diagram and $c$ is the upper reduction of $tb$, then $a\lambda c$ is an oriented circle diagram and
	\[
	\deg(a\lambda t\mu b) = \deg(a\lambda c)+\cu(t) + p-q,
	\]
	where $p$ (resp.~$q$) is the number of upper circles that are clockwise (resp.~anticlockwise) in the diagram $a\lambda t\mu b$.
\end{lem}
\begin{proof}
	When passing from $tb$ to $c$, we remove all upper circles, which obviously have the same number of cups and caps, and upper lines, which have one more cap than cup.
	From every other component we remove an equal number of cups and caps.
	The total number of cups removed is $\cu(t)$.
	The statement then follows from \cref{onemoreoneless}. 
\end{proof}
\begin{defi}
	Let $\Lambda$ and $\Gamma$ be two blocks, and let $t$ be a $\Lambda\Gamma$-matching.
	Define $K^{t}_{\Lambda\Gamma}$ to be the graded vector space with homogeneous basis
	\begin{equation*}
		\{(a\lambda t\mu  b)\mid\text{ for all closed oriented $\Lambda\Gamma$-circle diagrams } a\lambda t\mu  b\}.
	\end{equation*}
	Define a degree preserving linear map
	\begin{equation}\label{definvolution}
		*\colon K^{t}_{\Lambda\Gamma}\to K^{t ^*}_{\Gamma\Lambda}, \quad (a\lambda t\mu  b)\mapsto (b^*\mu t^*\lambda a^*),
	\end{equation}
	where $t^*$, $a^*$ and $b^*$ denote the mirror images of $t$, $a$, $b$ in the horizontal axis.
	
	We define $(c\nu d)\cdot(a\lambda t\mu b)$ to be $0$ if $d^*\neq a$.
	Otherwise, we draw $c\nu d$ beneath $a\lambda t\mu b$ and use the same surgery procedures as for $\KD$ to smooth out the symmetric middle section.
	The product is then defined as this linear combination.
	Similarly, we can define $(a\lambda t\mu b)\cdot(c\nu d)$ and thus endow $K^t_{\Lambda\Gamma}$ with the structure of a $(\KB_\Lambda, \KB_\Gamma)$-bimodule. 
\end{defi}

The following theorem generalizes the cellular structure of Khovanov's algebra to our setting and is very important for the following computations.
\begin{thm}\label{cellularstructuregeombimod}Suppose that we are given basis vectors $(a\lambda b)\in \KD_\Lambda$ and $(c\mu t\nu d)\in K^{t}_{\Lambda\Gamma}$.
	The multiplication satisfies
	\begin{equation}
		({a}\lambda b)(c\mu t\nu{d}) = \begin{cases}
			0&\text{if } b\neq c^*,\\
			s_{{a}\lambda b}(\mu)({a}\mu t\nu d)+(\dagger)&\text{if } b=c^*\text{ and }{a}\mu\text{ is oriented},\\
			(\dagger)&\text{otherwise,}
		\end{cases}
	\end{equation}
	where
	\begin{enumerate}
		\item $(\dagger)$ denotes a linear combination of basis vectors of $K^{t}_{\Lambda\Gamma}$ of the form $(a\nu't\nu'd)$ for $\mu'>\mu$,
		\item the scalar $s_{{a}\lambda b}(\mu)\in\{0,1,-1\}$ depends only on ${a}\lambda b$ and $\mu$, but not on ${d}$ and is equal to the scalar $s_{a\lambda b}(\mu)$ from \cite{ES2}*{Theorem 7.1} and
		\item if $\lambda=\mu$ and $b=\overline{\lambda}=c^*$, then $s_{a\lambda b}(\mu)=1$.
	\end{enumerate}
\end{thm}
\begin{proof}
	This follows by using the same arguments as \cite{ES2}*{Theorem 7.1 + Remark 7.6}.
\end{proof}

\subsection{Projective functors}
In this section we are going to introduce so-called projective functors and use them to prove that $\KD_\lambda$ is Koszul.
In \cite{BS3} the Khovanov algebra of type $\A$ is identified with a projective generator of parabolic category $\mathcal{O}$ for $\lie{gl}(n)$.
Under this identification the projective functors correspond to translation functors (see also \cite{HNS24} for the setting in type $\B$).

\begin{defi}\label{defiprojfunctor}
	Let $t$ be a proper $\Lambda\Gamma$-matching.
	Define the functor
	\begin{equation*}
		G^t_{\Lambda\Gamma}\coloneqq K^t_{\Lambda\Gamma}\langle-\ca(t)\rangle\otimes\_\colon \gMod_{lf}(\KB_\Gamma)\to \gMod_{lf}(\KB_\Lambda).
	\end{equation*}
	We call any functor which is isomorphic to a finite direct sum of the above functors (possibly shifted) a \emph{projective functor}.
\end{defi}
\begin{thm}\label{projfunctorsonproj}
	Let $t$ be a proper $\Lambda\Gamma$-matching and let $\gamma\in\Gamma$.
	Then
	\begin{enumerate}[label=(\arabic*)]
		\item\label{projfunctorsonproji} $G^t_{\Lambda\Gamma}P(\gamma)\cong K^t_{\Lambda\Gamma}e_\gamma\langle-\ca(t)\rangle$ as left $\KB_\Lambda$-modules,
		\item\label{projfunctorsonprojnonzero} the module $G^t_{\Lambda\Gamma}P(\gamma)$ is nonzero if and only if each upper line in $t\gamma\overline{\gamma}$ is oriented and
		\item\label{projfunctorsonprojexplicit} in this case moreover,
		\begin{equation*}
			G^t_{\Lambda\Gamma}P(\gamma)\cong P(\lambda)\otimes R^{\otimes n}\langle\cu(t)-\ca(t)\rangle
		\end{equation*} as graded left $\KB_\Lambda$-modules ($\KB_\Lambda$ acts again on the right-hand side only on the first factor), where $\lambda\in\Lambda$ is such that $\overline{\lambda}$ is the upper reduction of $t\overline{\gamma}$ and $n$ denotes the number of upper circles removed in the reduction process.  
	\end{enumerate}
\end{thm}
\begin{proof}
	For \cref{projfunctorsonproji} note that 
	\begin{align*}
		G^t_{\Lambda\Gamma}P(\gamma)&=K^t_{\Lambda\Gamma}\langle-\ca(t)\rangle\otimes_{\KB_\Gamma} P(\gamma) = K^t_{\Lambda\Gamma}\otimes_{\KB_\Gamma} \KB_\Gamma e_\gamma\langle-\ca(t)\rangle\\&\cong K^t_{\Lambda\Gamma}e_\gamma\langle-\ca(t)\rangle.
	\end{align*}
	For the forward implication of \cref{projfunctorsonprojnonzero}, note that for any weight $\nu$ such that $\nu\overline{\gamma}$ is oriented, the rays are oriented in the same ways as in $\gamma\overline{\gamma}$.
	Thus, if there exists an upper line in $t\gamma\overline{\gamma}$ which is not oriented, then there cannot exist an oriented $\Lambda\Gamma$-circle diagram of the form $a\mu t\nu\overline{\gamma}$.
	But these form a basis of $K^t_{\Lambda\Gamma}e_{\gamma}$ and hence $G^t_{\Lambda\Gamma}P(\gamma)=0$ by \cref{projfunctorsonproji}.
	
	In order to finish the proof, suppose that each upper line of $t\gamma\overline{\gamma}$ is oriented properly.
	Enumerate the $n$ upper circles in some order and define the map
	\begin{equation*}
		f\colon K^t_{\Lambda\Gamma}e_{\gamma}\to \KB_\Lambda e_\lambda\otimes R^{\otimes n}, \quad (a\mu t\nu \overline{\gamma})\mapsto (a\mu\overline{\lambda})\otimes x_i\otimes\dots\otimes x_n
	\end{equation*}
	where $x_i$ is $1$ (resp.~$X$) if the $i$-th circle is oriented anticlockwise (resp.~clockwise). 
	
	This map is then an isomorphism of vector spaces.  It is $\KB_\Lambda$-linear as every tag gets altered by an even number of undotted arcs, and
	moreover it is homogeneous of degree $\cu(t)$ by \cref{degreeformulaupperreduction}.
	By observing that $0\neq P(\lambda)=\KB_\Lambda e_\lambda$ and $G^t_{\Lambda\Gamma}P(\gamma)\cong K^t_{\Lambda\Gamma}e_\gamma\langle-\ca(t)\rangle$, this finishes the proof of \cref{projfunctorsonprojnonzero} and \cref{projfunctorsonprojexplicit}.
\end{proof}
\begin{cor}\label{geombimodareproj}
	The module $K^t_{\Lambda\Gamma}$ is projective as a left $\KB_\Lambda$-module as well as projective as a right $\KB_\Gamma$-module.
	Hence, $G^t_{\Lambda\Gamma}$ is exact.
\end{cor}
\begin{proof}
	By \cref{projfunctorsonproj}\cref{projfunctorsonproji} and \cref{projfunctorsonprojexplicit} we have that $K^t_{\Lambda\Gamma} = \bigoplus_{\gamma\in\Gamma}K^t_{\Lambda\Gamma}e_\gamma$ is projective as a left $\KB_\Lambda$-module.
	Using the antimultiplicative map $*$, $K^t_{\Lambda\Gamma}$ being projective as a right $\KB_\Gamma$-modules is the same as $K^{t^*}_{\Gamma\Lambda}$ being a projective left $\KB_\Gamma$-module, but this was done above.
\end{proof}

The following theorem deals with the effect of a projective functor on standard modules $\Delta(\mu)$.
\begin{thm}\label{projfunctorsoncell} The $\KB_\Lambda$-module $G^t_{\Lambda\Gamma}\Delta(\gamma)$ has a filtration
	\begin{equation*}
		\{0\}=M(0)\subset M(1)\subset\dots\subset M(n)=G^t_{\Gamma\Lambda}\Delta(\gamma)
	\end{equation*}
	such that $M(i)/M(i-1)\cong \Delta(\mu_i)\langle\deg(\mu_it\gamma)-\ca(t)\rangle$ for each $i$.
	In this case $\mu_1, \dots, \mu_n$ denote the elements of the set $\{\mu\in\Lambda\mid \mu t\gamma\text{ oriented}\}$ ordered such that $\mu_i>\mu_j$ implies $i<j$.
\end{thm}
\begin{proof}
	A basis for $G^t_{\Lambda\Gamma}P(\gamma)$ is given by $(a\mu t\nu\overline{\gamma})\otimes e_\gamma$ by \cref{projfunctorsonproj}.
	The module $G^t_{\Lambda\Gamma}U(\gamma)$ has a basis given by $(a\mu t\nu\overline{\gamma})\otimes e_\gamma$ with $\nu>\gamma$ by the paragraph after \cref{defkhovanovalg} and \cref{cellularstructuregeombimod}.
	Thus, a basis for $\Delta(\gamma)$ is given by $(a\mu t\gamma\overline{\gamma})\otimes e_\gamma$ 
	
	Now define $M(0)=\{0\}$ and inductively $M(i)$ to be the subspace generated by $M(i-1)$ and $\{(a\mu_it\gamma\overline{\gamma})\otimes e_\gamma\mid\text{ for all oriented cup diagrams } a\mu_i\}$.
	This defines a filtration of $G^t_{\Lambda\Gamma}\Delta(\gamma)$ by vector spaces with $M(n)=G^t_{\Lambda\Gamma}\Delta(\gamma)$ by the above argumentation.
	That the $M(i)$ are in fact $\KB_\Lambda$-submodules follows from \cref{cellularstructuregeombimod}, our assumption on the ordering of the $\mu_i$, and the above paragraph.
	
	The quotient $M(i)/M(i-1)$ has a basis given by
	\begin{equation*}
		\{(c\mu_it\gamma\overline{\gamma}\otimes e_\gamma\mid\text{ for all oriented cup diagrams } c\mu_i\}.
	\end{equation*}
	\cref{cellularstructuregeombimod} says that 
	\begin{equation*}
		(a\lambda b)(c\mu_it\gamma\overline{\gamma})\otimes e_\gamma\equiv\begin{cases}
			s_{a\lambda b}(\mu_i)(a\mu_it\gamma\overline{\gamma})\otimes e_\gamma&\text{if $b=c^*$ and $a\mu_i$ oriented,}\\
			0&\text{otherwise,}
		\end{cases}
	\end{equation*}
	working modulo $M(i-1)$.
	Looking at the definition of $\Delta(\mu)$, we see that the map
	\begin{equation*}
		M(i)/M(i-1)\to \Delta(\mu_i),\quad (c\mu_it\gamma\overline{\gamma})\otimes e_\gamma\mapsto(c\mu_i\overline{\mu_i})
	\end{equation*} is an isomorphism of $\KB_\Lambda$-modules.
	Moreover, it is homogeneous of degree $\deg(\mu_it\gamma)-\ca(t)$ by definition.
	This proves the first statement.
\end{proof}

\begin{thm}\label{standardkoszul}
	For every standard module $\Delta(\lambda)$ there exists a linear projective resolution
	\begin{equation*}
		\dots\to P^k\to P^{k-1}\to\dots\to P^1\to P^0\to \Delta(\lambda).
	\end{equation*}
\end{thm}
\begin{proof}
	This proof follows the same arguments as in the type $\A$ case in \cite{BS2}*{Theorem 5.3}.
	
	The claim is shown by a nested induction.
	First we do an induction on the number of caps of $\overline{\lambda}$ and secondly one on the Bruhat order.
	If $\#(\ca(\overline{\lambda}))=0$, then $\lambda$ is maximal.
	Hence, it suffices to consider maximal weights for the induction beginning.
	But in this case we have $P(\lambda)=\Delta(\lambda)$ and the claim holds. 
	
	Now suppose that $\#(\ca(\overline{\lambda}))>0$ and assume the claim for all $\mu$ with fewer caps and all $\lambda'>\lambda$.
	Observe that the number of caps is finite for every weight diagram.
	
	As $\lambda$ is not maximal, we can apply a Bruhat move $B$ at positions $i$ and $j$ to $\lambda$.
	This corresponds to a cap $C$ in $\overline{\lambda}$.
	Now let $\mu$ be the weight, which is obtained by deleting the positions of the endpoints of $C$, and denote the corresponding block by $\Gamma$.
	
	Then let $t$ be the $\Lambda\Gamma$-matching given by a cap connecting positions $i$ and $j$ (with the same parity of dots as $C$) and vertical strands everywhere else.
	
	We observe that there are exactly two weights $\gamma$ such that $\gamma t\mu$ is oriented, namely $\gamma=\lambda$ and $\gamma=\lambda'$ (the former corresponds to orienting $C$ anticlockwise and the latter to a clockwise orientation).
	So by \cref{projfunctorsoncell} (as $\lambda'>\lambda$) there is a short exact sequence
	\begin{equation}\label{sesstandardkoszul}
		\begin{tikzcd}
			0\arrow[r]&\Delta(\lambda')\arrow[r, "f"]&G^t_{\Lambda\Gamma}\Delta(\mu)\arrow[r]&\Delta(\lambda)\langle -1\rangle\arrow[r]&0.
		\end{tikzcd}
	\end{equation}
	By the induction hypothesis we have constructed a linear projective resolution $P^\bullet(\lambda')$ of $\Delta(\lambda')$ and as $\overline{\mu}$ contains fewer caps than $\overline{\lambda}$ we also have constructed a linear projective resolution $P^\bullet(\mu)$ of $\Delta(\mu)$.

	By \cref{geombimodareproj}, applying $G^t_{\Lambda\Gamma}$ to $P^\bullet(\mu)$ gives a projective resolution $G^t_{\Lambda\Gamma}P^\bullet(\mu)$ of $G^t_{\Lambda\Gamma}\Delta(\mu)$.
	
	As in \cite{BS2}*{Theorem 5.3}, the cone of $f$ is then a projective resolution of $\Delta(\lambda)\langle -1\rangle$, see also \cite{TS}*{Proposition 4.3.1} for the setting here.
	Now note that 
	\begin{equation*}
		\mathrm{Cone}(f)^{k+1}=P(\lambda')^{k}\oplus G^t_{\Lambda\Gamma}P(\mu)^{k+1}.
	\end{equation*}
	Now recalling the degree shift in \eqref{sesstandardkoszul}, we only need to show that $P(\lambda')^{k}$ and $G^t_{\Lambda\Gamma}P(\mu)^{k+1}$ are generated in degree $k$.
	For $P(\lambda')^k$, this follows as $P^\bullet(\lambda')$ is a linear projective resolution of $\Delta(\lambda')$.
	Furthermore, $P(\mu)^{k+1}$ is generated in degree $k+1$ and so $G^t_{\Lambda\Gamma}P(\mu)^{k+1}$ is generated in degree $k$ by \cref{projfunctorsonproj}.
\end{proof}
\begin{thm}\label{koszulb}
	The algebra $\KB$ is Koszul.
\end{thm}
\begin{proof}
	The category $\gMod_{lf}(\KB)$ is an upper-finite highest weight category by \cref{upperfinitehighweight}.
	A linear projective resolution of the left standard modules was constructed in \cref{standardkoszul}.
	As we have a duality preserving irreducibles, the statement follows from \cref{standardKoszulwithduality,koszulcrit}.
\end{proof}
\begin{rem}
	The theory of projective functors and geometric bimodules was further studied in \cite{HNS24}.
	In particular, there are analogues of \cref{projfunctorsonproj} for the effect on irreducible modules.
	The relevance of this formula is that they describe
	\begin{enumerate}
		\item the structure of indecomposable summands of $V^{\otimes d}$, where $V$ is the natural representation of $\Osp[r][2n]$, and
		\item the effect of translation functors on irreducible $\Osp[r][2n]$-modules.
	\end{enumerate}
\end{rem}

\subsection{Deligne categories}

In \cite{D19} Deligne  constructed families of universal tensor categories $\Rep(\mathrm{GL}_{\delta})$, $\Rep(\Or_{\delta})$, $\Rep(S_{\delta})$, $\delta \in \mathbb{C}$, \cite{D19} interpolating the representation categories $\Rep(GL(n))$, $\Rep(O(n))$, $\Rep(S(n))$, $n \in \mathbb{N}$. The interpolation property refers to the fact that these Deligne categories admit for integral parameters $\delta \in \mathbb{N}$ or $\mathbb{Z}$ the usual representation categories as quotients. We focus now on the orthogonal case introduced in \cite{D19} and further studied in \cite{CH}. A very similar situation exists in the $\Rep(\mathrm{GL}_{\delta})$-case.

The Deligne category $\Rep(\Or_{\delta})$ is generated by one object $V$ as a Karoubian symmetric monoidal category. The endomorphism spaces of tensor powers $V^{\otimes d}$ of this object are given by the Brauer algebras $\Br_d(\delta)$.


Let $\Lambda$ denote the set of all partitions, and  \begin{align*} \Lambda_d = \begin{cases}  \{ \lambda \ | \ |\lambda| = d - 2i, \ 0 \leq i \leq \frac{d}{2} \}, \text{ if } \delta \neq 0 \text{ or } d \text{ odd;} \\  \{ \lambda  \ | \ |\lambda| = d - 2i, \ 0 \leq i < \frac{d}{2} \}, \text{ if } \delta = 0 \text{ and } d \text{ even.} \end{cases} \end{align*} For any $\delta \in \mathbb{C}$, the algebras $Br_d(\delta)$ are cellular \cite{GL}. Its cell modules are denoted by $\Delta_{d,\delta}(\lambda)$ where $\lambda \in \Lambda_d(\delta) \cup \{ \emptyset \}$ and the indecomposable projective modules by $P_{d,\delta}(\lambda)$.

Using this parametrization of the indecomposable projective modules it is easy to show that the indecomposable objects in $\Rep(\Or_{\delta})$ are (up to isomorphism) parametrized by the set $\Lambda$ of all partitions. We denote by $R_{\delta}(\lambda)$ the corresponding indecomposable object. It arises as the image of an idempotent in some Brauer algebra $\Br_d(\delta) = \End(V^{\otimes d})$. If $R_{\delta}(\lambda) = Im(e_{\lambda}) \subset V^{\otimes d}$, $R_{\delta}(\mu) = Im(e_{\mu}) \subset V^{\otimes d}$, we have \[ \Hom_{\Rep(\Or_{\delta})}(R_{\delta}(\lambda),R_{\delta}(\mu)) = e_{\lambda} \Br_d(\delta) e_{\mu} \eqqcolon e_{\lambda}Ae_{\mu}.\]   

A representation of $\Rep(\Or_{\delta})$ is a contravariant functor from $\Rep(\Or_{\delta})$ to $Vect$, the category of finite dimensional complex vector spaces. We denote this category by $\mathcal{D}(\delta)$. This is an abelian category whose indecomposable projective objects are given by the representations \[ P_{\delta}(\lambda) \coloneqq \Hom_{\Rep(\Or_{\delta})} (-,R_{\delta}(\lambda))\]

The category $\mathcal{D}(\delta)$ can be identified with the category of locally finite dimensional $A$-modules for the locally finite dimensional locally unital algebra \[ A \coloneqq\bigoplus_{\lambda,\mu} e_{\lambda} A e_{\mu}.\] Under this identification, $P_{\delta}(\lambda)$ corresponds to the $A$-module $Ae_{\lambda}$. 

The category $\mathcal{D}(\delta)$ contains for each partition $\lambda$ representations $\Delta_{\delta}(\lambda)$ isomorphic to the representation which sends $R_{\delta}(\mu)$ to $\Hom_{Br_d(\delta)} (P_{d,\delta}(\mu),\Delta_{d,\delta}(\lambda))$ for some $d$ satisfying $\lambda,\mu \in \Lambda_d(\delta)$, and a morphism $\alpha\colon R_{\delta}(\nu) \to R_{\delta}(\nu)$ to the precomposition $-\circ \alpha$ with $\alpha$. By \cite[Corollary 2.11]{OSPII} the category $\mathcal{D}(\delta)$ is an upper finite highest weight category with poset $\Lambda$ and with standard objects $\Delta_{\delta}(\lambda), \lambda \in \Lambda$. Here we have fixed the reverse inclusion ordering on $\Lambda$: $\lambda \geq \mu$ if $\lambda$ is contained in $\mu$ (so that $\emptyset$ is maximal).

\begin{thm} \label{thm:del-and-k} \cite{OSPII} There is an equivalence of categories between $\mathcal{D}(\delta)$ and $\gMod_{lf}(\KB_{\Lambda})$. 
\end{thm}

This equivalence can be used to endow $\mathcal{D}(\delta)$ with a grading. The same arguments work if the Khovanov algebra of type $B$ is replaced with its type $A$ variant, and the Deligne category $\Rep(\Or_{\delta})$ with $\Rep(\mathrm{GL}_{\delta})$. The following statement follows from \cref{koszulb}.

\begin{cor} The graded version of $\mathcal{D}(\delta)$ is Koszul, i.e. $A = \bigoplus_{\lambda,\mu \in \Lambda} e_{\lambda} A e_{\mu}$ is a Koszul algebra. Likewise, the graded version of the abelianized Deligne category of $\Rep(\mathrm{GL}_{\delta})$ is Koszul.
\end{cor}

\begin{rem}
While the Khovanov algebras of type $A$ and $B$ admit a Koszul grading, there are also negative examples. In \cite{N23} Nehme constructed an analogue of the Khovanov algebra $K_n$ for representations of the Lie superalgebra $\mathfrak{p}(n)$ and used this to give an explicit description of the endomorphism ring of a projective generator for $\mathfrak{p}(n)$. By \cite[Theorem E]{N23} there does not exist a non-negative grading on $K_n$ with semisimple degree $0$ part, that is generated in degree $1$ for $n \geq 2$. In particular, $K_n$ cannot be Koszul. 
\end{rem}

\bibliography{biblio}
\end{document}